\documentclass[10pt,draft]{amsart}
\usepackage{amsmath,amssymb,amsrefs}

\numberwithin{equation}{section}

\newtheorem{lemma}[equation]{Lemma}

\newtheorem{prop}[equation]{Proposition}

\theoremstyle{definition}
\newtheorem{definition}[equation]{Definition}

\theoremstyle{remark}
\newtheorem{remark}[equation]{Remark}

\begin{document}

\begin{center}
\texttt{Comments, suggestions, corrections, and references are
  most welcomed!}
\end{center}

\title{On the orders of generators of capable $\mathbf{p}$-groups}
\author{Arturo Magidin}
\address{Dept. of Mathematical Sciences, The University of
  Montana, Missoula MT 59812}
\email{magidin@member.ams.org}

\subjclass[2000]{Primary 20D15}

\begin{abstract}
A group is called capable if it is a central factor group. For
each prime $p$ and positive integer~$c$, we prove the existence of a capable
$p$-group of class~$c$ minimally generated by an element of order $p$
and an element of order $p^{1+\lfloor\frac{c-1}{p-1}\rfloor}$. 
This is best possible.
\end{abstract}

\maketitle

\section{Introduction}\label{sec:intro}

Recall that a group $G$ is said to be \textbf{capable} if and only if
$G$ is isomorphic to $K/Z(K)$ for some group~$K$, where $Z(K)$ is the
center of~$K$. There are groups
which are not capable (nontrivial cyclic groups being a well-known
example), so capability places restrictions on the structure
of a group; see for example~\cites{heinnikolova,isaacs}. As noted by
P.~Hall in his landmark paper on the classification of
$p$-groups~\cite{hallpgroups}, the question of which $p$-groups
are capable is interesting and plays an important role in their
classification.

P.~Hall observed that if $G$ is a capable $p$-group of class~$c$, with
$c<p$, and $\{x_1,\ldots,x_n\}$ is a minimal set of generators with
$o(x_1)\leq o(x_2)\leq\cdots\leq o(x_n)$ (where $o(g)$ denotes the
order of the element~$g$), then $n>1$ and $o(x_{n-1})=o(x_n)$. 

In~\cite{capable} we used commutator calculus to derive a similar
necessary condition after dropping the hypothesis $c<p$: if $G$ is a
capable $p$-group of class~$c>0$, minimally generated by
$\{x_1,\ldots,x_n\}$, where $o(x_1)\leq\cdots\leq o(x_n)$, then we
must have $n>1$ and letting $o(x_{n-1}) = p^a$ and
$o(x_n)=p^b$, then $a$ and $b$ must satisfy
\begin{equation}
b \leq a + \left\lfloor\frac{c-1}{p-1}\right\rfloor,
\label{eq:inequality}
\end{equation}
where $\lfloor x\rfloor$ is the greatest integer less than or equal
to~$x$ (Theorem~3.19 in~\cite{capable}).  The dihedral group of order
$2^{c+1}$ shows that $(\ref{eq:inequality})$ is best possible when
$p=2$. The purpose of this note is to show that the inequality is 
best possible for all primes~$p$, thus answering in the affirmative
Question~3.22 in~\cite{capable}.

Notation will be standard; all groups will be written
multiplicatively, and we will denote the identity by~$e$. We use the
convention that the commutator of two elements $x$ and~$y$ is
$[x,y]=x^{-1}y^{-1}xy$. The lower central series of~$G$ is defined
recursively by letting $G_1=G$, and $G_{n+1}=[G_n,G]$; we say $G$ is
nilpotent of class (at most)~$c$ if and only if $G_{c+1}=\{e\}$. 
It is well known that if $G$ is of
class exactly~$c$, then $G_c\subset Z(G)$, and $G/Z(G)$ is nilpotent
of class exactly~$c-1$.

We let $C_n$ denote the cyclic group of order~$n$, and
$\mathbf{Z}$ the infinite cyclic group, both written
multiplicatively. 

\section{The case $c=1+(r-1)(p-1)$}\label{sec:construction}

The construction in this section is based on the example
given by Easterfield in Section~4 of~\cite{easterfield}.

Let $p$ be a prime, $r$ a positive integer. We construct a
$p$-group $K$ of class $c+1=2+(r-1)(p-1)$, minimally generated by an
element $y$ of order~$p$, and an element $x_0$ of order $p^{r}$. We
will show that the images of $y$ and $x_0$ have the same order in
$K/Z(K)$, thus exhibiting a capable group of class $c=1+(r-1)(p-1)$,
minimally generated by an element of order $p$ and one of order
$p^{1+\lfloor\frac{c-1}{p-1}\rfloor}$.

Let $H$ be the abelian group
\[ H = C_{p^r}\times C_{p^r} \times
\underbrace{C_{p^{r-1}}\times\cdots\times C_{p^{r-1}}}_{\mbox{$p-2$
    factors}}.\] 
Denote the generators of the cyclic factors of~$H$ by
$x_0,x_1,\ldots,x_{p-1}$, respectively. If $r=1$, then
    $x_1,\ldots,x_{p-1}$ are trivial. Let $y$ generate a cyclic
    group of order $p$, and let $y$ act on $H$ by $y^{-1}x_i y = x_i
    x_{i+1}$ for $0\leq i\leq p-2$ (so $[x_i,y]=x_{i+1}$), and
\[ y^{-1}x_{p-1}y = x_1^{-\binom{p}{1}}x_2^{-\binom{p}{2}}\cdots
    x_{p-2}^{-\binom{p}{p-2}} x_{p-1}^{1-\binom{p}{p-1}};\] as usual,
$\binom{n}{k}$ is the binomial coefficient $n$ choose~$k$. 
Let $K=H\rtimes \langle y\rangle$.

\begin{remark} The group constructed by Easterfield is the subgroup of
  $K$ generated by $y$ and $x_1,\ldots,x_r$. We can also realize $K$
  as the semidirect product of this subgroup by $\langle x_0\rangle$,
  letting $x_0$ act on $y$ by $x_0^{-1}yx_0 = yx_1^{-1}$, and act
  trivially on the~$x_i$.
\end{remark}

Note that $K$ is metabelian of class exactly $2+(r-1)(p-1)$. To verify
the class, note that $[K,K] = \langle x_1,\ldots,x_{p-1}\rangle$. We
then have:
\begin{eqnarray*}
K_3 & = & \langle x_1^p,x_2,\ldots,x_{p-1}\rangle;\\
K_4 & = & \langle x_1^p,x_2^p,x_3,\ldots,x_{p-1}\rangle;\\
&\vdots&\\
K_{2+(p-1)} & = & \langle x_1^p,x_2^p,\ldots,x_{p-1}^p\rangle;\\
K_{2+(p-1) + 1} & = & \langle x_1^{p^2},x_2^p,\ldots,x_{p-1}^p\rangle;\\
&\vdots&\\
K_{2+k(p-1)} & = & \langle
x_1^{p^k},x_2^{p^k},\ldots,x_{p-1}^{p^k}\rangle;\\
&\vdots&\\
K_{2+(r-1)(p-1)} & = & \langle
  x_1^{p^{r-1}},x_2^{p^{r-1}},\ldots,x_{p-1}^{p^{r-1}}\rangle =
  \langle x_1^{p^{r-1}}\rangle.
\end{eqnarray*}
Finally, note that $x_1^{p^{r-1}}$ is central: $y^{-1}x_1^{p^{r-1}}y =
(x_1x_2)^{p^{r-1}}=x_1^{p^{r-1}}$. Therefore $K$ is of class exactly $2+(r-1)(p-1)$.

The group $G=K/Z(K)$ will therefore be of class $1+(r-1)(p-1)$,
minimally generated by $yZ(K)$ and $x_0Z(K)$. The order of $yZ(K)$ is
of course equal to~$p$. As for $x_0Z(K)$, note that no nontrivial
power of $x_0$ is central: if $x_0^k$ is central, then
\[ x_0^k = y^{-1}x_0^ky = (y^{-1}x_0y)^k = (x_0x_1)^k = x_0^kx_1^k;\]
therefore $x_1^k = e$, which implies that $p^r|k$,
so $x_0^k = e$. Therefore, the order of $x_0Z(K)$ is
$p^{r}$. Thus, $G$ is a capable group of class $c$, with $c=1+(r-1)(p-1)$,
minimally generated by an element of order $p$ and an element of order
$p^{r} = p^{1+\lfloor\frac{c-1}{p-1}\rfloor}$.

We note the following fact about $K$, which we will use in
the following section:

\begin{lemma} Let $p$ be any prime, and let $r$ be an arbitrary
  positive integer. There exists a group $K$ of class $2+(r-1)(p-1)$,
  generated by elements $y$ and $x_0$ of orders $p$ and $p^r$, respectively,
  such that $x_0^{p^{r-1}}$ does not commute with $y$.
\label{lemma:notcentral}
\end{lemma}

\section{General case}\label{sec:induction}

Again, let $p$ be a prime, and let $c>1$ be an arbitrary
integer. We want to exhibit a capable group~$G$ of class exactly~$c$,
generated by an element of order $p$ and an element of order
$p^{1+\lfloor\frac{c-1}{p-1}\rfloor}$. 

Our construction in this section will be based on the
nilpotent product of groups; we specialize the definition to the case
we are interested in:

\begin{definition} Let $A_1,\ldots,A_n$ be cyclic groups, and let
  $c>0$. The $c$-nilpotent product of the~$A_i$, denoted
  $A_1\amalg^{{\germ N}_c}\cdots\amalg^{{\germ N}_c} A_n$ is defined
  to be the group $F/F_{c+1}$, where $F$ is the free product of the
  $A_i$, $F = A_1 * \cdots * A_n$, and $F_{c+1}$ is the $(c+1)$-st
  term of the lower central series of~$F$.
\end{definition}

It is easy to verify that the $c$-nilpotent product of the $A_i$ is of
class exactly~$c$, and that it is
their coproduct (in the sense of category theory) in the variety
${\germ N}_c$ of all nilpotent groups of class at most~$c$. The
$1$-nilpotent product is simply the direct sum of the~$A_i$.

Note that if $G$ is the $c$-nilpotent
product of the $A_i$, then $G/G_{k+1}$ is the $k$-nilpotent product
of the $A_i$ for all $k$, $1\leq k\leq c$.

We consider $\mathcal{G}=C_p\amalg^{{\germ N}_{c+1}}\mathbf{Z}$, the
$(c+1)$-nilpotent product of a cyclic group of order $p$ and the
infinite cyclic group. Denote the generator of the finite cyclic group
by $a$, and the generator of the infinite cyclic group by $z$. Let
$G=\mathcal{G}/Z(\mathcal{G})$. Then $G$ is capable of class~$c$. We
want to show that $zZ(\mathcal{G})$ has the required order.

\begin{prop} Let $a$ generate $C_p$ and $z$ generate the infinite
  cyclic group~$\mathbf{Z}$. If $\mathcal{G}=C_p\amalg^{{\germ N}_{c+1}}
  \mathbf{Z}$, then 
\[ Z(\mathcal{G})\cap \langle z\rangle = \Bigl\langle
  z^{p^{1+\lfloor\frac{c-1}{p-1}\rfloor}}\Bigr\rangle.\]
\end{prop}

\begin{proof} The fact that $z^{p^{1+\lfloor\frac{c-1}{p-1}\rfloor}}$ is central
  follows from Theorem~3.16 in~\cite{capable}, so we just need to
  prove the other inclusion.  We proceed by induction on~$c$.  The
  claim is true if $c=1$ since the commutator bracket is bilinear in a
  group of class two. Assume
  the inclusion holds for $c-1$, with $c>1$. Note that $\langle
  z\rangle\cap \mathcal{G}_2 = \{e\}$.

Consider $\mathcal{G}/\mathcal{G}_{c+1}$; this is the $c$-nilpotent
product of $C_p$ and $\mathbf{Z}$, so by the induction hypothesis, the
intersection of the center and the subgroup generated by $z$ is
generated by the $p^{1+\lfloor\frac{c-2}{p-1}\rfloor}$-st power
of~$z$. Since the center of $\mathcal{G}$ is contained in the pullback
of the center of $\mathcal{G}/\mathcal{G}_{c+1}$, we deduce that the
smallest power of $z$ that could possibly be in $Z(G)$ is the
$p^{1+\lfloor\frac{c-2}{p-1}\rfloor}$-st power.

If $\lfloor\frac{c-2}{p-1}\rfloor = \lfloor\frac{c-1}{p-1}\rfloor$,
then we are done. So the only case that needs to be dealt with is the
case considered in the previous section, when $c = 1 + (r-1)(p-1)$ for
some positive integer~$r>1$.

Here we use the universal property of the coproduct. Let $K$ be the
group from Lemma~\ref{lemma:notcentral}. Since $\mathcal{G}$ is the
coproduct of $C_p$ and $\mathbf{Z}$ in ${\germ N}_{c+1}$, the morphisms
$C_p\to K$ given by $a\mapsto y$, and $\mathbf{Z}\to K$
given by $z\mapsto x_0\in K$, induce a unique
homomorphism $\varphi\colon\mathcal{G}\to K$. The image of
$Z(\mathcal{G})$ must lie in $Z(K)$ (since the map is surjective).
Since $\varphi(z^{p^{r-1}}) =
x_0^{p^{r-1}}$ does not commute with $y$, we conclude that
$z^{p^{r-1}}\not\in Z(\mathcal{G})$. This proves that the smallest
power of~$z$ that could lie in $Z(\mathcal{G})$ is $z^{p^{r}}$,
 which gives the desired inclusion.
\end{proof}

Now let $G = \mathcal{G}/Z(\mathcal{G})$. This is a group of
class~$c$, minimally generated by $aZ(\mathcal{G})$ and
$zZ(\mathcal{G})$. The former has order~$p$, and the latter element has order
$p^{1+\lfloor\frac{c-1}{p-1}\rfloor}$ by the proposition above. Thus $G$ is a capable group
of class~$c$, minimally generated by two elements whose orders satisfy
the equality in~(\ref{eq:inequality}), showing that the inequality is
indeed best possible.

\begin{remark} I believe that in general
  inequality~$(\ref{eq:inequality})$ will be both necessary and sufficient
  for the capability of a $c$-nilpotent product of cyclic $p$-groups. This
  is indeed the case when $c<p$ and when $p=c=2$; see~\cite{capable}. However,
I have not been able to establish this
for arbitrary $p$ and $c$, which forced the somewhat indirect approach
taken in this note.
\end{remark}

\section*{Acknowledgements}

I thank Prof. Avinoam Mann for bringing Easterfield's
paper~\cite{easterfield} to my attention. It proved invaluable in
suggesting the construction in Section~\ref{sec:construction}.


\section*{References}

\begin{biblist}
\bib{easterfield}{article}{
  author={Easterfield, T.E.},
  title={The orders of products and commutators in prime-power groups},
  date={1940},
  journal={Proc. Cambridge Philos. Soc.},
  volume={36},
  pages={14\ndash 26},
  review={\MR {1,104b}},
}
\bib{hallpgroups}{article}{
  author={Hall, P.},
  title={The classification of prime-power groups},
  date={1940},
  journal={J. Reine Angew. Math},
  volume={182},
  pages={130\ndash 141},
  review={\MR {2,211b}},
}
\bib{heinnikolova}{article}{
  author={Heineken, Hermann},
  author={Nikolova, Daniela},
  title={Class two nilpotent capable groups},
  date={1996},
  journal={Bull. Austral. Math. Soc.},
  volume={54},
  number={2},
  pages={347\ndash 352},
  review={\MR {97m:20043}},
}
\bib{isaacs}{article}{
  author={Isaacs, I. M.},
  title={Derived subgroups and centers of capable groups},
  date={2001},
  journal={Proc. Amer. Math. Soc.},
  volume={129},
  number={10},
  pages={2853\ndash 2859},
  review={\MR {2002c:20035}},
}
\bib{capable}{article}{
  author={Magidin, Arturo},
  title={Capability of nilpotent products of cyclic groups},
  eprint={{arXiv:math.GR/0403188}},
  note={Submitted},
}

\end{biblist}
\end{document}